\documentclass{llncs}
\usepackage{epsfig}
\usepackage{amssymb}
\usepackage{epic, eepic}
\usepackage{makeidx}
\usepackage{mathrsfs} 
\usepackage{braket}
\usepackage{cite} 

\newcommand{\bl}{\mbox{b}}
\newcommand{\bn}{\mbox{bn}}
\newcommand{\pn}{\mbox{pn}}
\newcommand{\fpn}{\mbox{fpn}}
\newcommand{\minb}{\mbox{minb}}

\newtheorem{thm}{Theorem}
\newtheorem{lem}[thm]{Lemma}
\newtheorem{prop}[thm]{Proposition}
\newtheorem{cor}[thm]{Corollary}
\newtheorem{ex}[thm]{Example}
\newtheorem{df}[thm]{Definition}
\newtheorem{con}[thm]{Conjecture}

\newcommand{\tume}{\hspace*{0.2em}}

\def\NZQ{\mathbb}               

\def\RR{{\NZQ R}}

\def\SCR{\mathscr}
\def\sC{{\SCR C}}

\newenvironment{pf}{\begin{quotation}%
        \par\hspace{-10mm}\makebox[10mm][l]{Proof:}}{%
        \hspace*{\fill} \hbox{\rule{6pt}{6pt}}\end{quotation}\par}
\author{Kenji Kashiwabara \inst{1} \and Tadashi Sakuma\inst{2}}
\institute{Department of General Systems Studies \\ University of Tokyo \and Systems Science and Information Studies \\%
Faculty of Education, Art and Science \\%
Yamagata University \\%
}
\title{On ideal minimally non-packing clutters}
\begin{document}
\maketitle

\begin{abstract}
We consider the following conjecture proposed by Cornu\'ejols, Guenin and Margot: every ideal minimally non-packing clutter has a transversal of size 2. For a clutter ${\cal C}$, let $\tilde{\cal C}$ denote the set of hyperedges of $\cal C$ which intersect any minimum transversal in exactly one element. We divide the (non-)existence problem of an ideal minimally non-packing clutter $\cal D$ into two steps. In the first step, we give necessary conditions for ${\cal C}=\tilde{\cal D}$ when ${\cal D}$ is an ideal minimally non-packing clutter. In the second step, for a clutter $\cal C$ satisfying the conditions in the first step, we consider whether $\cal C$ has an ideal minimally non-packing clutter $\cal D$ with ${\cal C}=\tilde{\cal D}$.
We show that the clutter of a combinatorial affine plane satisfies the conditions in the first step. Moreover, we show that the clutter of a combinatorial affine plane does not have any ideal minimally non-packing clutter of blocking number at least 3.
\end{abstract}

\section{Introduction}

\subsection{Background and motivation}
In the celebrated paper \cite{Seymour_1977} of Seymour, motivated by 
the pluperfect and (weak) perfect graph theorems for the set covering 
problem by Fulkerson and Lov\'{a}sz, he introduced the concept 
of so-called ``the Max-Flow-Min-Cut property'' of clutters, which 
is the packing counterpart of the totally dual integrality built in the 
perfection. That is, a clutter $\sC$ has the {\em Max-Flow-Min-Cut 
property} (the MFMC property, for short) if, for its clutter matrix 
$M(\sC)$, the linear system $M(\sC) \vec{x} \geqq  \vec{1}, \vec{x} 
\geqq \vec{0}$ is totally dual integral. 
(A matrix inequality $A \vec{x} \geq \vec{b}$ (resp. to $A \vec{x} \leq 
\vec{b}$) is called {\it totally dual integral} if the linear program 
$\min\{\langle \vec{w}, \vec{x}\rangle | A \vec{x} \geq \vec{b}\}$
(resp. to $\max\{\langle \vec{w}, \vec{x}\rangle | A \vec{x} \leq
\vec{b}\}$) has an integral optimal dual solution $\vec{y}$ for every 
integral cost vector $\vec{w}$ for which the above linear program has 
a finite optimum.)
In the case of the anti-blocking polytope of a clutter matrix, 
its integrality and the totally dual integrality of its linear system are 
coincident with the perfection. Seymour also pointed out that this 
``obvious analog'' of the set covering problem is false for the set 
packing problem, because there exists a non-MFMC clutter 
$Q_6:=\big\{\{1,3,5\},\{1,4,6\},\{2,3,6\},\{2,4,5\}\big\}$ whose 
blocking polyhedron 
$\set{\vec{x}\in \RR^{6} |  \vec{0}\leq \vec{x}, M({Q_6})\vec{x}\geq \vec{1}}$ 
is integral (i.e. {\em ideal\/}). On the other hand, he proved that 
this $Q_6$ is the only ideal binary clutter which is minimally non-MFMC 
as the meaning of clutter minor. 

A clutter $\sC$ has the {\em packing property\/} (resp. {\em packs\/})
if the both sides of the linear programming
equation $\min\{ \langle \vec{\omega}, \vec{x} \rangle | 
\vec{x} \geqq \vec{1}, M(\sC) \vec{x} \geqq \vec{1}\} 
= \max\{ \langle \vec{y},  \vec{1} \rangle  | \vec{y} \geqq \vec{0}, 
\vec{y} M(\sC) \leqq \vec{\omega}\}$ have optimal solution integral 
vectors $\vec{x}$ and $\vec{y}$ for all cost vectors $\vec{\omega}$ 
with components equal to $0,1$ or $\infty$ (resp. when 
$\vec{\omega} = \vec{1}\/$). Lehman \cite{Lehman_1979} proved that 
the packing property implies the idealness. However, the converse is 
false because the ideal clutter $Q_6$ does not pack and hence does 
not have the packing property. By definition, the MFMC property 
implies the packing property. But how about the converse? In 1993, 
Conforti and Cornu\'{e}jols \cite{Conforti_and_Cornuejols_1993} 
proposed the following famous conjecture. 
\begin{con}[Conforti and Cornu\'{e}jols 1993]\label{C&C}
{\rm A clutter has the packing property if and only if it has the MFMC 
property.}
\end{con}
Despite its natural appearance, this conjecture is very difficult 
and still open. The existing approaches can be classified roughly 
into two categories: 

\noindent
The first category is to find a clutter class for which 
the conjecture is affirmative (or, if possible, false). 
Conjecture\tume\ref{C&C} holds for the binary clutters 
\cite{Seymour_1977}, the diadic clutters 
\cite{Cornuejols_Guenin_and_Margot_2000}, 
the clutter of circuits of a digraph \cite{Guenin_2001}, 
the clutter of cycles in an undirected graph \cite{Ding_and_Zang_2002}, 
the broken circuit clutter of two-dimension affine 
convex geometries \cite{Hachimori_and_Nakamura_2008}, 
the Ehrhart clutters \cite{Martinez-Bernal_O'Shea_and_Villarreal_2010}, 
and so on. However, for the almost all of them, except for the 
case of the binary clutters shown in the Seymour's initial paper 
\cite{Seymour_1977}, the MFMC property is coincident with not only 
the packing property but also the idealness. In other words, there 
are only minimally non-ideal excluded clutter minors for these 
classes to have the MFMC property. Of course, there are several 
known clutter classes, other than the binary clutters, on which 
the packing property is not coincident with the idealness. (It is 
well known that the clutter of dijoins \cite{Schrijver_1980, 
Cornuejols_and_Guenin_2002, Williams_and_Guenin_2005} falls 
into the case. See also \cite{Kashiwabara_and_Sakuma_2010} 
for another example.)
However, as far as the authors know, Conjecture\tume\ref{C&C} 
seems unsettled even if restricted to each of these classes. To begin 
with, there are so few clutter classes on which the packing 
property is characterized by the set of minimally excluded minors 
which inevitably includes some ideal clutters (again, see  
\cite{Kashiwabara_and_Sakuma_2010}). 

\noindent
The second category is to investigate ``the packing property'' itself 
and extract key nature of the concept by which we can 
prove or disprove the conjecture. The first essential 
step on this line was achieved by Cornu\'{e}jols, Guenin 
and Margot \cite{Cornuejols_Guenin_and_Margot_2000}. 
Starting with the discovery of $Q_6$ \cite{Seymour_1977}, 
there have been found numerous (and several infinite families of) ideal 
minimally non-packing clutters until today (e.g. \cite{Schrijver_1980, 
Guenin_1998, Cornuejols_Guenin_and_Margot_2000, 
Cornuejols_and_Guenin_2002, Williams_and_Guenin_2005, 
Kashiwabara_and_Sakuma_2010}). All of these existing clutters have 
the common property: The blocking number is $2$ for all of them. 
Cornu\'{e}jols, Guenin and Margot
\cite{Cornuejols_Guenin_and_Margot_2000} 
conjectured that the converse is also true. 
\begin{con}\label{C&G&M}
The blocking number of every ideal minimally non-packing clutter 
is $2$. 
\end{con}
\vspace*{-10pt}
\noindent
Furthermore, they proved that the above conjecture implies 
Conjecture\tume\ref{C&C}. 

In this paper, the authors will provide a framework to attack 
Conjecture\tume\ref{C&G&M}. A {\em tilde core} 
of an ideal minimally non-packing clutter $\sC$ is the maximal set of 
hyperedges of $\sC$ such that every minimum transversal of $\sC$ 
has a unique common element with each of the hyperedges. 
As the concept of the {\em core} has greatly developed the theory 
of minimally non-ideal clutters (see \cite{Cornuejols_2001} for details), 
the concept of the tilde core may have similar impact to the theory 
of ideal minimally non-packing clutters. Actually, Cornu\'{e}jols, 
Guenin and Margot \cite{Cornuejols_Guenin_and_Margot_2000} proved 
that several key features of the ideal minimally non-packing clutters 
are controlled by their tilde cores. The authors will develop their idea 
to a framework to check whether a given clutter can be a tilde core 
of an ideal minimally non-packing clutter or not. This framework is 
useful not only for the search of counterexamples but also to prove 
the conjecture. We demonstrate this by applying our framework to 
the case of a special clutter, namely, the combinatorial affine planes.  
We show that every combinatorial affine plane whose blocking number is at least 3 
 cannot be a tilde core of any ideal minimally non-packing clutter (Theorem \ref{thm:affinenon}).
In connection with this, we should note that whether 
a combinatorial projective plane except for the Fano plane 
$F_7$ can be a core of a minimally non-ideal clutter or not 
is a famous open question of the theory of minimally non-ideal 
clutters (see Question\tume$6$ in 
\cite{Cornuejols_Guenin_and_Tuncel_2008}).

\subsection{Overview of our results} 

We consider Conjecture \ref{C&G&M} in this paper. That is, we consider the (non-)existence problem of an ideal minimally non-packing clutter of blocking number at least 3. We propose a new framework to attack the conjecture. 


Let $E$ be a finite ground set of clutters throughout this paper.
$\tilde{\cal C}$ denotes the set of hyperedges in a clutter $\cal C$ each of which intersects any minimum transversal in exactly one element. 
A tilde clutter $\tilde{\cal C}$ was first introduced in Cornu\'ejols, Guenin and Margot \cite{Cornuejols_Guenin_and_Margot_2000}. That paper gave necessary conditions for $\cal C$ to be an ideal minimally non-packing clutter in terms of $\tilde{\cal C}$. In our paper, we develop their idea. We contrive tractable necessary conditions for $\cal C$ to be an ideal minimally non-packing clutter in terms of $\tilde{\cal C}$. By our approach, clutters that we have to consider are restricted.

We divide the (non-)existence problem of an ideal minimally non-packing clutter $\cal D$ as in Conjecture \ref{C&G&M} into two steps. In the first step (Section \ref{sec:tilde}), we give necessary conditions for ${\cal C}=\tilde{\cal D}$ when ${\cal D}$ is an ideal minimally non-packing clutter. 
We call a clutter satisfying the conditions in the step 1 a precore clutter. In the second step (Section \ref{sec:solution}), for a precore clutter $\cal C$, we consider whether $\cal C$ has an ideal minimally non-packing clutter $\cal D$ with ${\cal C}=\tilde{\cal D}$. 
Since the necessary conditions in step 1 are rather strong, clutters that we have to consider are much confined.
However, we found a several classes of precore clutters. When we try to find a counterexample or prove the conjecture, we have only to consider the problem for each precore clutter $\cal C$. That is, it is the problem for $\cal C$ to have an ideal non-minimally non-packing clutter $\cal D$ with ${\cal C}=\tilde{\cal D}$. Starting with a (rather vague) task to find a counterexample to Conjecture \ref{C&G&M}, here, we obtain a tractable concrete problem whether a given clutter $\cal C$ has an ideal non-minimally non-packing clutter $\cal D$ with ${\cal C}=\tilde{\cal D}$ or not.

Section \ref{sec:tilde} is devoted to step 1. For an ideal non-packing clutter $\cal D$, we present several necessary conditions of $\tilde{\cal D}$: the integral blocking condition, tilde-invariance, the integrality of $\mbox{I}(\tilde{\cal D})$, and non-separability (Theorem \ref{thm:step1}). 
The integral blocking condition is defined as the coincidence of the fractional packing number and the blocking number. 
This condition is a fundamental condition as a premise of an argument.
A clutter $\cal C$ satisfying the integral blocking condition is called tilde-invariant if ${\cal C}=\tilde{\cal C}$ holds.
For an ideal clutter $\cal C$, $\tilde{\cal C}$ is tilde-invariant.
A clutter is ideal if and only if the blocking polyhedron $\{x\in {\RR}^E|\langle 1_H,x\rangle \geq 1 \mbox{ for all }H\in {\cal C},x\geq 0\}$ is an integral polyhedron.
The polyhedron $\mbox{I}({\cal C})$ is a face of the above blocking polyhedron defined by the equalities corresponding to minimum transversals. We show that, for an ideal clutter $\cal C$, not only $\mbox{I}({\cal C})$ but also $\mbox{I}(\tilde{\cal C})$ is an integral polyhedron (Theorem \ref{thm:icint2}).
The minimum transversals define the affine hull of $\mbox{I}({\cal C})$ and some non-minimum transversals define facets of $\mbox{I}({\cal C})$. By observing these transversals carefully, we can derive useful information from them.

Cornu\'ejols, Guenin and Margot \cite{Cornuejols_Guenin_and_Margot_2000} proved that deleting all the elements on a hyperedge of an ideal minimally non-packing clutter decreases the blocking number by at least two. We call such a condition hyperedge-nonseparability.
We also present a condition called non-separability, which is a generalization of hyperedge-separability.

We have the following implications among conditions on $\cal C$ under the conditions that its minimum transversals cover $E$ and the integral blocking condition (Lemma \ref{lem:justone}, Lemma \ref{lem:itos} and Lemma \ref{lem:dim}).

Integrality of $\mbox{I}({\cal C})$ $\Rightarrow$ tilde-full condition+dimension condition $\Rightarrow$ \\
tilde-full condition $\Rightarrow$ weak tilde-invariant clutter.


In Section \ref{sec:solution}, when a precore clutter $\cal C$ is given, we present several necessary conditions for an ideal minimally non-packing clutter ${\cal D}$ with ${\cal C}=\tilde{\cal D}$: Conditions IM, IF, H, and B (Theorems \ref{thm:step2} and \ref{thm:step2b}). Since these conditions for $\cal D$ are strong enough, the next condition for a precore $\cal C$ is derived. When a precore clutter $\cal C$ has an ideal minimally non-packing clutter, there must exist a clutter $\cal D$ satisfying Conditions IM, IF, H, and B  (Corollary \ref{cor:step2}).  

In Section \ref{sec:unique}, we consider the problem with an additional condition that the maximum fractional packing is unique.
Many classes of precore clutters satisfy this condition as far as we know.
In this case, $I(\tilde{\cal C})$ for a precore clutter is an integral simplex. This condition is characterized in terms of transversals and a condition about dimension (Theorem \ref{thm:simplexfacet}).
We give an example of a precore clutter, namely, a combinatorial affine plane. The clutter ${\cal C}$ of a combinatorial affine plane is obtained by deleting one element from a combinatorial projective plane. 
We show that the clutter $\cal C$ of a combinatorial affine plane is a precore clutter (Theorem \ref{thm:cppintegral}). Moreover, we show that the clutter $\cal C$ of a combinatorial affine plane cannot have a counterexample $\cal D$ to Conjecture \ref{C&G&M} with ${\cal C}=\tilde{\cal D}$ (Theorem \ref{thm:affinenon}).

\section{Preliminaries}

Let $E$ be a finite set.
A family ${\cal C}\subseteq 2^E$ of sets is said to be a {\it clutter} if no member includes another member. A member of $\cal C$ is called a {\it hyperedge}. For details about clutters, please refer to \cite{Cornuejols_2001}.

For a clutter $\cal C$, a set on $E$ is a {\it transversal} if it intersects every element of $\cal C$ and it is minimal with respect to inclusion in such sets. 
$\bl({\cal C})$ denotes the clutter consisting of all the transversals of ${\cal C}$.
A {\it minimum transversal} of $\cal C$ is a transversal of the minimum size.
$\minb({\cal C})$ denotes the set of minimum transversals of a clutter ${\cal C}$.
Note that we assume that the word `transversal' always means a `minimal' transversal to avoid the confusion between a minimum transversal and a minimal transversal in our definition.

The {\it blocking number} $\bn({\cal C})$ of a clutter ${\cal C}$ is the minimum size of a transversal in $\bl({\cal C})$. 
The {\it packing number} $\pn({\cal C})$ of a clutter ${\cal C}$ is the maximum size of a family of hyperedges of a clutter such that any pair of them does not intersect. Clearly, 
$\pn({\cal C})\leq \bn({\cal C})$ holds. When $\pn({\cal C})=\bn({\cal C})$ holds, $\cal C$ is said to {\it pack}. When ${\cal C}\subseteq {\cal C}'$, $\pn({\cal C})\leq \pn({\cal C}')$ and $\bn({\cal C})\leq \bn({\cal C}')$ hold.

The {\it contraction} of $A$ from $\cal C$ is ${\cal C}/A=\min(\{X-A|X\in {\cal C}\})$ where min is the operation of collecting minimal sets with respect to inclusion.
The {\it deletion} of $A$ from $\cal C$ is ${\cal C}\backslash A=\{X\in {\cal C}|X\cap A=\emptyset \}$. A {\it minor} of $\cal C$ is a clutter which is obtained by contractions and deletions iteratively from $\cal C$. A {\it proper minor} means a minor which is not equal to the original clutter. The {\it restriction} of $\cal C$ to $A$ is ${\cal C}[A]={\cal C}\backslash A^c$.

A clutter $\cal C$ is called  {\it minimally non-packing} if it does not pack and any proper minor packs.
A clutter $\cal C$ is called  {\it minimally non-packing with respect to deletion} if it does not pack and any proper deletion minor packs.
A clutter on $E$ is called {\it minimum-transversal-covered} if its minimum transversals cover $E$.

\begin{lem}\label{lem:cover}
For any minimally non-packing clutter with respect to deletion, it is minimum-transversal-covered.
\end{lem}

\begin{pf}
Since $\cal C$ does not pack, $\pn({\cal C}) < \bn({\cal C})$ holds. For a minimally non-packing clutter $\cal C$  with respect to deletion and $a\in E$, the deletion ${\cal C}\backslash a$ packs. Therefore $\pn({\cal C}\backslash a) = \bn({\cal C}\backslash a)$. Since deleting one element decreases the blocking number by at most one, we have $\bn({\cal C}\backslash a) = \bn({\cal C})-1$. Recall that the deletion of a clutter corresponds to the contraction of the clutter of its transversals. If $a$ is not covered by any minimum transversal, every minimum transversal of ${\cal C}\backslash a$ is also a minimum transversal of $\cal C$, a contradiction to $\bn({\cal C}\backslash a) = \bn({\cal C})-1$.
\end{pf}

For a clutter $\cal C$, $M({\cal C})$ denotes a clutter matrix of $\cal C$, whose row vectors coincide with the incidence vectors of its hyperedges. We consider the following linear problem. $$\max\{\sum_{H\in {\cal C}}y(H) |yM({\cal C})\leq 1_{\cal C},y\in {\RR}^{\cal C}\}= \min\{\sum_{a\in E}x(a)|M({\cal C})x \geq 1_E,x\in {\\R}^E\}.$$ Note that the above equality always holds because of the duality theorem of the linear programming. We call the maximizing problem of $y$ the {\it primal problem} and the minimizing problem of $x$ the {\it dual problem}. 

A clutter ${\cal C}$ is {\it ideal} if $\{x\in {\RR}^E|\langle x, 1_H \rangle \geq 1 \mbox{ for all }H\in {\cal C}, x \geq 0\}$ is an integral polyhedron where $1_H$ is the incidence vector of $H$. Note that $x \geq 0$ means $x(a)\geq 0$ for any $a\in E$. It is known that any minor of an ideal clutter is an ideal clutter again.

By the complementary slackness of the linear programming, for a maximum solution $y$ of the primal problem, a minimum solution $x$ of the dual problem and any $a\in E$, $x(a)>0$ implies $\sum_{H:a\in H\in {\cal C}}y(H) = 1$.

A {\it maximum fractional packing} $y$ of a clutter $\cal C$ is a function ${\cal C}\to {\RR}_\geq$ maximizing the sum $\sum_{H\in {\cal C}}y(H)$ such that $\sum_{H:a\in H\in {\cal C}}y(H)\leq 1$ for any $a\in E$. Every maximum fractional packing is an optimal solution of the primal problem. The {\it support} of a maximum fractional packing $y$ is the set of hyperedges $H$ with $y(H)>0$. Define
$$\mbox{F}({\cal C})=\{z\in \RR^{{\cal C}}|z \mbox{ is a maximum fractional packing}\}.$$

The {\it fractional packing number} is $\sum_{H\in{\cal C}} y(H)$ for $y\in \mbox{F}({\cal C})$, denoted by $\fpn({\cal C})$. Note that $\bn({\cal C})\geq \fpn({\cal C}) \geq \pn({\cal C})$. For an ideal clutter $\cal C$, $\fpn({\cal C})=\bn({\cal C})$ holds. When a clutter ${\cal C}
$ packs, $\pn({\cal C})=\fpn({\cal C})=\bn({\cal C})$.

$\tilde{\cal C}$ denotes the set of hyperedges in a clutter $\cal C$ which intersect any minimum transversal in exactly one element. That is,
$$\tilde{\cal C}=\{H\in {\cal C}:|B\cap H|=1 \mbox{ for all }B\in \minb({\cal C})\}.$$
We call $\tilde{\cal C}$ the {\it tilde clutter} of ${\cal C}$.
A clutter $\tilde{\cal C}$, obtained by the tilde operation, plays a crucial role in this paper.

\section{Precore conditions}\label{sec:tilde}

In this section, we present several necessary conditions for $\tilde{\cal D}$  when $\cal D$ is ideal minimally non-packing: the integral blocking condition, the integrality of $\mbox{I}({\cal C})$, and non-separability. 

\subsection{Integral blocking condition}

\begin{df}\normalfont
A clutter ${\cal C}$ satisfies the {\it integral blocking condition} if its fractional packing number $\fpn({\cal C})$ is equal to its blocking number $\bn({\cal C})$.
\end{df}

\begin{lem}\label{lem:slack}
Assume that a minimum-transversal-covered clutter ${\cal C}$ satisfies the integral blocking condition. Then, for any $y\in \mbox{F}({\cal C})$, $\sum_{H\in {\cal C}}y(H)1_H=1_E$ holds.
\end{lem}

\begin{pf}
By the integral blocking condition, we have $\sum_{H\in {{\cal C}}}y(H)=\bn({\cal C})$ for a maximum fractional packing $y$. By the complementary slackness, we have $\sum_{H\in {{\cal C}}}y(H)1_H=1_E$ since the minimum transversals cover $E$. Note that its minimum transversals are optimal solutions of the dual problem by the integral blocking condition.
\end{pf}

\begin{lem}\label{lem:ibc}
Assume that a minimum-transversal-covered clutter ${\cal C}$ satisfies the integral blocking condition. Then every hyperedge in the support of a maximum fractional packing of $\cal C$ intersects any minimum transversal in exactly one element. That is, any hyperedge in the support of some maximum fractional packing of $\cal C$ belongs to $\tilde {\cal C}$. 
\end{lem}

\begin{pf}
By definition of transversals, any hyperedge $H\in {\cal C}$ and any minimum transversal $B$ satisfy $|H\cap B| \geq 1$. When there exist some hyperedge $H$ and some minimum transversal $B$ with $|H\cap B|>1$, $\langle \sum_{H\in {{\cal C}}}y(H)1_H, 1_B\rangle=\sum_{H\in {{\cal C}}}y(H)\langle 1_H, 1_B\rangle > \sum_{H\in {{\cal C}}}y(H)$.
Since $1_E=\sum_{H\in {{\cal C}}}y(H)1_H$ by Lemma \ref{lem:slack}, $\langle 1_E, 1_B\rangle=\langle \sum_{H\in {{\cal C}}}y(H)1_H, 1_B\rangle>\sum_{H\in {{\cal C}}}y(H)=\bn({\cal C})$, which contradicts the fact that $\langle 1_E, 1_B\rangle=\bn({\cal C})$.
\end{pf}

So we can regard $y\in \mbox{F}({\cal C})$ as $y\in \mbox{F}(\tilde{\cal C})$.

\begin{lem}\label{lem:mfp}
Assume that a minimum-transversal-covered clutter $\cal C$ satisfies the integral blocking condition. 
For $y\in \mbox{F}({\cal C})$, $\sum_{H\in {\tilde{\cal C}}}y(H)1_H=1_E$ and $\sum_{H\in {\tilde{\cal C}}}y(H)=\bn(\tilde{\cal C})=\bn({\cal C})$. Moreover $\mbox{F}({\cal C})=\mbox{F}({\tilde{\cal C}})$ holds.
\end{lem}

\begin{pf}
By the integral blocking condition and covering by the minimum transversals, $\sum_{H\in {{\cal C}}}y(H)=\bn({\cal C})$ holds for $y\in \mbox{F}({\cal C})$. We have $\sum_{H\in {{\cal C}}}y(H)1_H=1_E$ by Lemma \ref{lem:slack}. 
Therefore  $\sum_{H\in {\tilde{\cal C}}}y(H)1_H=\sum_{H\in {{\cal C}}}y(H)1_H=1_E$ holds by Lemma \ref{lem:ibc}. By taking the inner product between each side of the equality and a minimum transversal $B$ of $\cal C$, we have $\sum_{H\in {\tilde{\cal C}}}y(H)\langle 1_H,1_B \rangle = \langle 1_E, 1_B \rangle$. Since $\langle 1_H,1_B\rangle =1$ for $H\in {\tilde{\cal C}}$, $\sum_{H\in {\tilde{\cal C}}}y(H)=\bn({\cal C})$. Since $\tilde{\cal C}\subseteq {\cal C}$, $\sum_{H\in {\tilde{\cal C}}}y(H)\leq\fpn(\tilde{\cal C})\leq \bn(\tilde{\cal C})\leq\bn({\cal C})$.
Therefore $y$ attains a maximum fractional packing of $\tilde{\cal C}$.
So $\mbox{F}({\cal C})=\mbox{F}({\tilde{\cal C}})$.
We have $\fpn(\tilde{\cal C})=\bn({\cal C})=\bn(\tilde{\cal C})$.
\end{pf}

\begin{cor}\label{cor:equal}
For a minimum-transversal-covered clutter $\cal C$ which satisfies the integral blocking condition, $\tilde{\cal C}$ is minimum-transversal-covered and also satisfies the integral blocking condition. Moreover, $\minb({\cal C})\subseteq \minb(\tilde{\cal C})$ holds.
\end{cor}

\begin{pf}
By Lemma \ref{lem:mfp}, $\mbox{F}({\cal C})=\mbox{F}({\tilde{\cal C}})$ and $\bn({\cal C})=\bn(\tilde{\cal C})$. So $\fpn(\tilde{\cal C})=\fpn({\cal C})=\bn({\cal C})=\bn(\tilde{\cal C})$.

By definition, $\tilde{\cal C}\subseteq \cal C$ holds. So every minimum transversal of $\cal C$ intersects every hyperedge of $\tilde{\cal C}$. Since $\bn({\cal C})=\bn(\tilde{\cal C})$, a minimum transversal of $\cal C$ is a minimum transversal of ${\tilde{\cal C}}$. So the minimum transversals of $\tilde{\cal C}$ cover $E$.
\end{pf}

\begin{ex}
Let ${\cal C}=\{ac,bc,bd\}$ on $E=\{a,b,c,d\}$. Then $\bl({\cal C})=\{ab,bc,cd\}$ and $\tilde{\cal C}=\{ac,bd\}$. Since $\bl(\tilde{\cal C})=\{ab,bc,cd,da\}$, this is an example whose minimum transversals are different between ${\cal C}$ and $\tilde{\cal C}$.
\end{ex}

\begin{cor}\label{cor:intblk}
Let $\cal C$ be an ideal minimum-transversal-covered clutter. Then $\cal C$ satisfies the integral blocking condition. Moreover $\tilde{\cal C}$ satisfies the integral blocking condition.
\end{cor}

\begin{pf}
By the duality theorem of linear programming, for any maximum fractional packing $y$ on $\cal C$, there exists a minimum transversal $B\in \minb({\cal C})$ with $yM({\cal C})1_B=|B|$. Note that we can take an integral optimal solution $1_B$ since $\cal C$ is ideal. 
By Lemma \ref{lem:ibc}, $yM({\cal C})1_B=yM(\tilde{\cal C})1_B$.
Since $M(\tilde{\cal C})1_B=1_{\tilde{\cal C}}$, $yM(\tilde{\cal C})1_B=\sum_{H\in \tilde{\cal C}}y(H)$ is also the fractional packing number of $\cal C$. 
Therefore $\cal C$ satisfies the integral blocking condition.

Moreover, by Corollary \ref{cor:equal}, $\tilde{\cal C}$ also satisfies the integral blocking condition.
\end{pf}




\begin{prop}
When every minor of a clutter satisfies the integral blocking condition, the clutter is an ideal clutter.
\end{prop}

\begin{pf}
When the clutter is not ideal, it has a minimally non-ideal clutter ${\cal C}'$ as a minor. Then $\fpn({\cal C}')> \bn({\cal C}')$ since any minimum transversal of ${\cal C}'$ is not any optimal solution of the dual problem.
The clutter ${\cal C}'$ does not satisfy the integral blocking condition.
 \end{pf}



\begin{lem}\label{lem:nonemp}
For a minimum-transversal-covered clutter $\cal C$ which satisfies the integral blocking condition, the number of hyperedges in $\tilde {\cal C}$ is at least the blocking number.
\end{lem}

\begin{pf}
There exists at least one maximum fractional packing. The number of hyperedges which belong to the support is at least the blocking number  since, for each minimum transversal $B$, every element of $B$ intersects a different hyperedge in $\tilde {\cal C}$. Therefore the statement follows from Lemma \ref{lem:ibc}. 
\end{pf}


\subsection{Tilde-invariant clutters and tilde-full condition}

\begin{df}\normalfont
A clutter $\cal C$ is a {\it tilde-invariant clutter} if it satisfies ${\cal C}=\tilde{\cal C}$ and the integral blocking condition.
\end{df}

\begin{df}\normalfont
A clutter $\cal C$ is  a {\it weak tilde-invariant clutter} if $\cal C$ satisfies the integral blocking condition and $\tilde{\cal C}$ is a tilde-invariant clutter.
\end{df}

By definition, every tilde-invariant clutter is a weak tilde-invariant clutter.

\begin{ex}\label{ex:kyou}
Even if the minimum transversals of a clutter cover $E$, and if it satisfies the  integral blocking condition, 
it may not be a weak tilde-invariant clutter.
Let ${\cal C}=\{abc,de,ef,fd,af,bd,ce\}$ on $\{a,b,c,d,e,f\}$ shown in Figure \ref{fig:nonintegral}. We have $\bl({\cal C})=\{ade,bef,cfd\}$, and $\tilde{\cal C}=\{af,bd,ce,abc\}$.
Since $\tilde{\cal C}$ is not a tilde-invariant clutter, $\cal C$ is not a weak tilde-invariant clutter.


\begin{figure}[ht]
\begin{center}
\includegraphics[width=4cm]{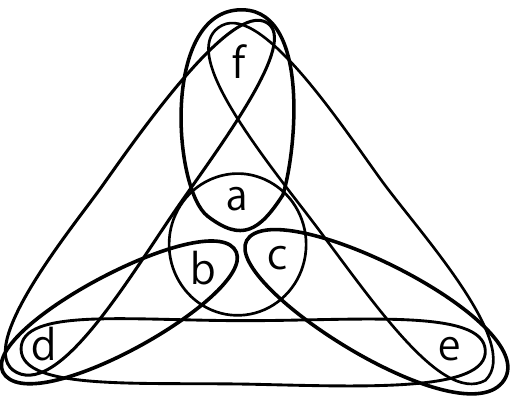}
\caption{an example of a non-weak tilde-invariant clutter}\label{fig:nonintegral}
\end{center}
\end{figure}

\end{ex}

\begin{df}\normalfont
A clutter $\cal C$ satisfies the {\it tilde-full condition} when $\cal C$ satisfies the following conditions.
\begin{itemize}
\item It is minimum-transversal-covered.
\item $\cal C$ satisfies the integral blocking condition.
\item Every hyperedge in $\tilde {\cal C}$ belongs to the support of some maximum fractional packing.
\end{itemize}
\end{df}




\begin{lem}\label{lem:justone}
If a clutter $\cal C$ satisfies the tilde-full condition, $\cal C$ is a weak tilde-invariant clutter.
\end{lem}

\begin{pf}
Assume that $\cal C$ satisfies the tilde-full condition. Then any hyperedge $H$ in $\tilde{\cal C}$ belongs to the support of some maximum fractional packing $y$. By Lemma \ref{lem:mfp}, $y$ is also a maximum fractional packing of $\tilde{\cal C}$. On $\tilde{\cal C}$, $H$ and any minimum solution $x_0\in\{x\in{\RR}^E|M(\tilde{\cal C})x\geq 1_E,x\geq 0\}$ in the dual problem  satisfy $\langle 1_H,x_0\rangle =1$. Since any minimum transversal $B\in \minb(\tilde{\cal C})$ is a minimum solution, we have $|B\cap H|=1$. Therefore $\tilde{\cal C}$ is a tilde-invariant clutter.
\end{pf}



\begin{lem}\label{lem:ibcibc}
When a clutter ${\cal C}$ satisfies the tilde-full condition, $\tilde {\cal C}$ also satisfies the tilde-full condition.
\end{lem}

\begin{pf}
By Corollary \ref{cor:equal}, $\tilde{\cal C}$ satisfies the integral blocking condition.
For any $H\in \tilde{\cal C}$, there exists a maximum fractional packing $y$ of $\cal C$ whose support contains $H$ since ${\cal C}$ satisfies the tilde-full condition. Since the support of $y$ is contained in  $\tilde{\cal C}$ by Lemma \ref{lem:mfp}, $y$ is also a maximum fractional packing of $\tilde{\cal C}$. Therefore $\tilde {\cal C}$  satisfies the tilde-full condition.
\end{pf}

\subsection{Polytope $\mbox{I}({\cal C})$}

We define a polyhedron $\mbox{I}({\cal C})$ as follows.
$$\mbox{I}({\cal C})=\{x\in {\RR}^E: \langle x,1_D \rangle=1 \mbox{ for all $D\in \minb$}({\cal C}), \langle x,1_D\rangle \geq 1 \mbox{ for all } D\in \bl({\cal C}), x \geq 0\}.$$

Note that this polyhedron $I({\cal C})$ is a face of the blocking polyhedron $\{x\in {\RR}^E: \langle x,1_D\rangle \geq 1 \mbox{ for all } D\in \bl({\cal C}), x \geq 0\}.$ This polyhedron plays a central role in this paper.

\begin{lem}\label{lem:nonempi}
For a clutter ${\cal C}$ which satisfies the integral blocking condition, $\mbox{I}({\cal C})$ is a non-empty polyhedron.
For a clutter ${\cal C}$ which satisfies the integral blocking condition, $\mbox{I}({\cal C})$ is a polytope if and only if $\cal C$ is minimum-transversal-covered. 
\end{lem}

\begin{pf}
We show the first statement.
Since $\cal C$ satisfies the integral blocking condition, there exists at least one maximum fractional packing whose support intersects any minimum transversal in exactly one element.
Therefore $\tilde C$ contains some hyperedge $H$. Then $1_H\in \mbox{I}({\cal C})$.

We show the second statement.

Assume that $\cal C$ is minimum-transversal-covered. 
By Lemma \ref{lem:nonemp}, $\tilde{\cal C}$ is non-empty. Since the incidence vector of every hyperedge of $\tilde{\cal C}$ belongs to $\mbox{I}({\cal C})$, $\mbox{I}({\cal C})$ is non-empty. 

Assume $x\in \mbox{I}({\cal C})$. For any $a\in E$, there exists $B\in \minb({\cal C})$ with $a\in B$ since the minimum transversals cover $E$. Therefore $x(a)$ is at most 1 since $x\geq 0$ and $\langle x,1_B\rangle=1$. Since $x \geq 0$, $\mbox{I}({\cal C})$ is bounded.

Conversely, assume that there exists $a\in E$ which is covered with no minimum transversals. Since $\cal C$ satisfies the integral blocking condition, $\tilde C$ contains some hyperedge $H$. Then $1_H\in \mbox{I}({\cal C})$. $1_H+k1_a \in \mbox{I}({\cal C})$ for any $k\geq 0$. So $\mbox{I}({\cal C})$ is not bounded.
\end{pf}

\begin{lem}\label{lem:icint}
For an ideal minimum-transversal-covered clutter ${\cal C}$, $\mbox{I}({\cal C})$ is an integral polytope.
\end{lem}

\begin{pf}
Since $\cal C$ is ideal, $\{x\in {\RR}^E|\langle x,1_B\rangle\geq 1 \mbox{ for any }B\in \bl({\cal C}),x\geq 0\}$ becomes an integral polyhedron. Since $\mbox{I}({\cal C})$ is a face of it, $\mbox{I}({\cal C})$ is also an integral polyhedron. Note that every face of an integral polyhedron is integral. $\mbox{I}({\cal C})$ is a polytope since its minimum transversals cover $E$ and Lemma \ref{lem:nonempi}.
\end{pf}



\begin{lem}\label{lem:tanten}
For a minimum-transversal-covered clutter ${\cal C}$ satisfying the integral blocking condition, all the integral points in $\mbox{I}({\cal C})$ coincide with the incidence vectors of $\tilde{\cal C}$.
\end{lem}

\begin{pf}
Every incidence vector of $\tilde{\cal C}$ satisfies all the inequalities defining $\mbox{I}({\cal C})$. Conversely, consider an integral point $x$ in $\mbox{I}({\cal C})$. Note that such an integral point $x$ is a 01-vector because the minimum transversals cover $E$. Let $H$ be the set with $1_H=x$. We show that $H$ is a hyperedge of $\tilde{\cal C}$. 
By the definition of $\mbox{I}({\cal C})$, $H$ intersects every transversal and intersects every minimum transversal in exactly one element.  Next we show that such $H$ is minimal.
If there exists another integral point $1_{H'}$ such that $H'\subsetneq H$. Then there exists $a\in H-H'.$ Since the minimum transversals cover $E$, there exists a minimum transversal $B$ containing $a$. But $B-a$ also intersects any hyperedge, which contradicts the minimality of a hyperedge. So such a set is a hyperedge in $\tilde{\cal C}$.
\end{pf}

\begin{lem}\label{lem:center}
For a minimum-transversal-covered clutter ${\cal C}$ which satisfies the integral blocking condition, the point consisting of all $1/\bn({\cal C})$ is contained in the relative interior of $\mbox{I}({\cal C})$.
\end{lem}

\begin{pf}
Since any transversal defining a facet of $\mbox{I}({\cal C})$ has a size of at least $\bn({\cal C})+1$, the inner product of a transversal defining a facet of $\mbox{I}({\cal C})$ and the point in the statement is more than 1. On the other hand, the inner product of a minimum transversal of $\mbox{I}({\cal C})$ and the point in the statement is exactly 1. Therefore the point consisting of $1/\bn({\cal C})$ is contained in the relative interior of $\mbox{I}({\cal C})$.
\end{pf}

\begin{lem}\label{lem:itos}
Assume that a clutter ${\cal C}$ which satisfies the integral blocking condition. When $\mbox{I}({\cal C})$ is an integral polytope, ${\cal C}$ satisfies the tilde-full condition.
\end{lem}

\begin{pf}
Since $\mbox{I}({\cal C})$ is a polytope, ${\cal C}$ is minimum-transversal-covered by Lemma \ref{lem:nonempi}.
By the integrality of the polytope $\mbox{I}({\cal C})$ and Lemma \ref{lem:tanten}, there exists no extreme point other than such incidence vectors of $\tilde{\cal C}$. Therefore the point in Lemma \ref{lem:center} is expressed as a positive combination of the extreme points of $\mbox{I}({\cal C})$ because the point is in the relative interior of $\mbox{I}({\cal C})$. By multiplying such a coefficient by $\bn({\cal C})$, the sum of all the components of the vector attains $\bn({\cal C})$, which is the fractional packing number. So there is a maximum fractional packing such that all the coefficients are positive. 
\end{pf}

\begin{thm}\label{thm:icint2}
For an ideal minimum-transversal-covered clutter ${\cal C}$, $\mbox{I}(\tilde{\cal C})$ is an integral polytope with $\mbox{I}({\cal C})=\mbox{I}(\tilde{\cal C})$. The extreme points of $\mbox{I}({\cal C})$ consist of the incidence vectors of $\tilde{\cal C}$.
\end{thm}

\begin{pf}
$\mbox{I}({\cal C})$ is an integral polytope by Lemma \ref{lem:icint}. By Corollary \ref{cor:equal} and Corollary \ref{cor:intblk}, $\cal C$ and $\tilde {\cal C}$ are minimum-transversal-covered and satisfy the integral blocking condition. Since $\tilde{\cal C}\subseteq {\cal C}$, there exists $B'\in \bl(\tilde{\cal C})$ with $B'\subseteq B$ for any $B\in \bl({\cal C})$. Therefore $I(\tilde{\cal C})\subseteq I({\cal C})$ holds. Therefore $\mbox{I}({\cal C})=\mbox{I}(\tilde{\cal C})$ follows from Lemma \ref{lem:tanten}.
\end{pf}

\begin{df}\normalfont
A clutter $\cal C$ satisfies the {\it dimension condition} if 
$$(\mbox{affine dimension of }\tilde{\cal C})+(\mbox{affine dimension of }\minb(\tilde{\cal C}))=|E|-1.$$
 The affine dimension of $\tilde{\cal C}$ means the dimension of the affine hull of all the incidence vectors of $\tilde{\cal C}$.
\end{df}

\begin{lem}\label{lem:dim}
For a tilde-invariant clutter ${\cal C}=\tilde{\cal C}$ such that $\mbox{I}({\cal C})$ is an integral polytope, ${\cal C}$ satisfies the dimension condition.
\end{lem}

\begin{pf}
Since $\mbox{I}({\cal C})$ is a polytope, ${\cal C}$ is minimum-transversal-covered by Lemma \ref{lem:nonempi}.
Since $\mbox{I}({\cal C})$ is an integral polytope, the extreme points of $\mbox{I}({\cal C})$ consist of the incidence vectors of $\tilde{\cal C}$ by Lemma \ref{lem:tanten}. 
So the dimension of $\mbox{I}({\cal C})$ is equal to the dimension of the affine hull  of $\tilde{\cal C}$. 
We have only to show that the dimension of $I({\cal C})$ is not affected by transversals other than $\minb(\tilde{\cal C)}$.
If the dimension of $I({\cal C})$ is affected by a non-minimum transversal, such a non-minimum transversal intersects any hyperedge in $\cal C$ in exactly one element. For a maximum fractional packing $y$, $\sum_{H\in \tilde{\cal C}}y(H)=\sum_{H\in \tilde{\cal C}}y(H)\langle 1_H,1_B \rangle = \sum_{H\in \tilde{\cal C}}y(H)\langle 1_H,1_B\rangle=\langle1_E,1_B\rangle=  |B| > \bn({\cal C})$, a contradiction to Lemma \ref{lem:mfp}.
\end{pf}


\subsection{Non-separability}

Separability is a necessary condition for a clutter $\cal C$ to have an ideal minimally non-packing clutter $\cal D$ with ${\cal C}=\tilde{\cal D}$.

\begin{df}\normalfont
A clutter ${\cal C}$ is {\it separable} if
there exists a non-trivial partition $\{E_1,E_2\}$ of $E$ and $\bn({\cal C})=\bn({\cal C}[E_1])+\bn({\cal C}[E_2])$ where ${\cal C}[E_1]:={\cal C}\backslash E_1^c$ with ${\cal C}[E_2]:={\cal C}\backslash E_2^c$. Otherwise it is called {\it non-separable}.
\end{df}

\begin{lem}\label{lem:nonsep}
A minimally non-packing clutter with respect to deletion is non-separable.
\end{lem}

\begin{pf}
Since a clutter $\cal C$ is minimally non-packing, $\cal C$ does not pack.
Assume that it is separable with a partition $\{E_1,E_2\}$ of $E$.

Consider the case where both ${\cal C}[E_1]$ and  ${\cal C}[E_2]$ pack.
Then $\fpn({\cal C}[E_1])=\pn({\cal C}[E_1])=\bn({\cal C}[E_1])$ and $\fpn({\cal C}[E_2])=\pn({\cal C}[E_2])=\bn({\cal C}[E_2])$ and $\fpn({\cal C})\geq \fpn({\cal C}[E_1])+\fpn({\cal C}[E_2])$. 
By separability, $\pn({\cal C})\geq \pn({\cal C}[E_1])+\pn({\cal C}[E_2])=\bn({\cal C}[E_1])+\bn({\cal C}[E_2])=\bn({\cal C})$. 
Since $\pn({\cal C})\leq \bn({\cal C})$ generally, $\pn({\cal C})=\bn({\cal C})$ holds. So $\cal C$ packs.
This contradicts the fact that $\cal C$ does not pack. 

Consider the case where either of them does not pack. This contradicts the assumption that $\cal C$ is minimally non-packing.
\end{pf}

\begin{lem}\label{lem:nonsep2}
Assume that a minimum-transversal-covered clutter $\cal C$ satisfies the integral blocking condition and non-separability.
Then $\tilde{\cal C}$ is non-separable.
\end{lem}

\begin{pf}
Assume that $\tilde{\cal C}$ is separable with a partition $\{E_1,E_2\}$ such that $\bn(\tilde{\cal C}[E_1])+\bn(\tilde{\cal C}[E_2])=\bn(\tilde{\cal C})$. We have $\bn(\tilde{\cal C}[E_1])\leq \bn({\cal C}[E_1])$ and $\bn(\tilde{\cal C}[E_2])\leq \bn({\cal C}[E_2])$ since $\tilde{\cal C}[E_1]\subseteq {\cal C}[E_1]$ and $\tilde{\cal C}[E_2]\subseteq {\cal C}[E_2]$. By Lemma \ref{lem:mfp} and the  integral blocking condition on $\cal C$, we have $\bn({\cal C})=\bn(\tilde{\cal C})$.  Since $\bn({\cal C}[E_1])+\bn({\cal C}[E_2])\leq \bn({\cal C})$ generally, 
we have $\bn(\tilde{\cal C})=\bn(\tilde{\cal C}[E_1])+\bn(\tilde{\cal C}[E_2])\leq\bn({\cal C}[E_1])+\bn({\cal C}[E_2])\leq \bn({\cal C})$. Therefore 
we have $\bn({\cal C}[E_1])+\bn({\cal C}[E_2])=\bn({\cal C})$, which contradicts the fact that $\cal C$ is non-separable.
\end{pf}



\begin{df}\normalfont
A minimum-transversal-covered clutter ${\cal C}$ satisfying the integral blocking condition is {\it hyperedge-separable} if there exists a hyperedge $H\in \tilde{\cal C}$ such that $\bn({\cal C}\backslash H)=\bn({\cal C})-1$. Otherwise, that is, when $\bn({\cal C}\backslash H)<\bn({\cal C})-1$ for any $H\in \tilde{\cal C}$, it is called {\it hyperedge-nonseparable}.
\end{df}

\begin{lem}
When a minimum-transversal-covered clutter ${\cal C}$ satisfying the integral blocking condition is hyperedge-separable, it is separable.
\end{lem}

\begin{pf}
When a clutter $\cal C$ is hyperedge-separable at $H\in \tilde{\cal C}$, we take a partition $\{H,H^c\}$ of $E$. Then we have $\bn({\cal C}[H])+\bn({\cal C}\backslash H)=\bn({\cal C})$ because of  $\bn({\cal C}[H])=1$.
\end{pf}

\subsection{Summarizing the conditions in step 1}\label{subsec:tilde}

\begin{df}\normalfont
For a clutter $\cal C$, a clutter $\cal D$ is called a {\it solution clutter} of $\cal C$ if ${\cal C}=\tilde{\cal D}$.
\end{df}

As a weaker problem, we firstly consider the problem for a clutter $\cal C$ to have an ideal clutter $\cal D$ with ${\cal C}=\tilde{\cal D}$. That is, we discard the condition of `minimally non-packing'.

We have considered necessary conditions for the tilde-invariant clutter $\cal C$ to have an ideal minimally non-packing clutter $\cal D$ with ${\cal C}=\tilde{\cal D}$. In this subsection, we integrate these results.

\begin{thm}\label{thm:step1}
Assume that a clutter $\cal C$ has an ideal minimally non-packing solution $\cal D$. Then $\cal C$ satisfies the following conditions.

\begin{itemize}
\item ${\cal C}$ satisfies the integral blocking condition.
\item $\mbox{I}({\cal C})$ is an integral polytope.
\item ${\cal C}$ is non-separable.
\end{itemize}
\end{thm}

\begin{pf}
The clutter ${\cal C}$ is minimum-transversal-covered by Lemma \ref{lem:cover}.
The clutter ${\cal C}$ is integral blocking by Corollary \ref{cor:intblk}.
So the integrality of $\mbox{I}({\cal C})$ follows from Theorem \ref{thm:icint2}. 
The non-separability condition follows from Lemmas \ref{lem:nonsep} and \ref{lem:nonsep2}.
\end{pf}

In other words, for any ideal minimally non-packing solution $\cal D$, $\tilde {\cal D}$ satisfies the conditions in Theorem \ref{thm:step1}.

We call a clutter ${\cal C}=\tilde{\cal C}$ satisfying the conditions in Theorem \ref{thm:step1} a {\it precore clutter}.
When we consider Conjecture \ref{C&G&M}, we have only to consider the precore clutters.
We give an example of a precore clutter in Subsection \ref{subsec:projective}.

\begin{thm}
Assume that ${\cal C}$ satisfies the integral blocking condition and $\mbox{I}({\cal C})$ is an integral polytope.
Then $\tilde{\cal C}$ is minimum-transversal-covered and tilde-invariant, and satisfies the tilde-full condition, and the dimension condition.
\end{thm}

\begin{pf}
The clutter $\tilde{\cal C}$ is minimum-transversal-covered by Lemma \ref{lem:nonempi}.
The integral blocking condition follows from Corollary \ref{cor:intblk}.
The clutter $\tilde{\cal C}$ satisfies the tilde-full condition by Lemma \ref{lem:itos} and Lemma \ref{lem:ibcibc}.
So $\tilde{\cal C}$ is a tilde-invariant clutter by Lemma \ref{lem:justone}.
The dimension condition follows from Lemma \ref{lem:dim}.
\end{pf}





\section{Conditions in the second step}\label{sec:solution}

After we find a clutter $\cal C$ satisfying the conditions in step 1, we have to discuss whether it has an ideal clutter $\cal D$ which is minimally non-packing with ${\cal C}=\tilde{\cal D}$ further.


For a precore clutter $\cal C$, we discuss necessary conditions for an ideal solution clutter $\cal D$.
The difference between the conditions in Theorem \ref{thm:step1} and those in this section is whether $\cal D$ appears in the conditions directly or not. 

{\bf Condition I:} $\mbox{I}({\cal C})=\mbox{I}({\cal D})$ holds.

We can divide Condition I into Conditions IM and IF.

{\bf Condition IM:} The affine space generated by the incidence vectors of the minimum transversals of $\cal C$ is equal to be the affine space generated by the incidence vectors of the minimum transversals of $\cal D$.

{\bf Condition IF:} If a facet $F$ of $\mbox{I}({\cal C})$ is defined by a transversal of ${\cal C}$,
there exists at least one element $B\in \bl({\cal D})$ defining the facet $F$.
Moreover, every $B\in \bl({\cal D})$ intersects any $H\in {\cal C}$.

{\bf Condition H:} ${\cal C} \subseteq{\cal D}$ must hold. For any $H\in {\cal D}-{\cal C}$, there exists $B\in \minb({\cal C})$ with $|H\cap B| \geq 2$.

\begin{thm}\label{thm:step2}
For a precore clutter $\cal C$, every ideal solution clutter $\cal D$ to $\cal C$ satisfies Conditions IM, IF, and H.
\end{thm}

\begin{pf}
By Theorem \ref{thm:icint2}, Condition I holds. Since the affine hull of $I({\cal C})$ is defined by $\minb({\cal C})$, Condition IM holds.
Since the facets of $I({\cal C})$ are defined by $\bl({\cal C})$ and $x\geq 0$, Condition IF holds.
Note that every $B\in \bl({\cal D})$ intersects any $H\in {\cal C}$ since ${\cal C}\subseteq {\cal D}$.


Assume $H\in {\cal D}-{\cal C}$. Since $H \notin \tilde{\cal C}={\cal C}$, there exists $B\in \minb({\cal C})$ with $|H\cap B| \geq 2$.
Therefore Condition H holds.


\end{pf}

We consider necessary conditions for the tilde-invariant clutter $\cal C$ to have an ideal minimally non-packing solution clutter $\cal D$.

{\bf Condition B:} For any disjoint sets $A,B \subseteq E$ with $A\cup B\neq \emptyset$, $\bn({\cal C}/A\backslash B)\leq\pn({\cal D}/A\backslash B)$ holds.

\begin{thm}\label{thm:step2b}
For a precore clutter $\cal C$, every ideal minimally non-packing solution clutter $\cal D$ to $\cal C$ satisfies Condition B.
\end{thm}

\begin{pf}
Since any proper minor of $\cal D$ has the packing property, $\bn({\cal D}/A\backslash B)=\pn({\cal D}/A\backslash B)$ holds. 
Since ${\cal C}\subseteq {\cal D}$, $\bn({\cal C}/A\backslash B)\leq \bn({\cal D}/A\backslash B)$. Therefore Condition B holds.
\end{pf}

\begin{cor}\label{cor:step2}
When a precore clutter $\cal C$ has an ideal minimally non-packing solution, there must exist a clutter $\cal D$ satisfying Conditions IF, IM, H, and B.
\end{cor}

We have not found a precore clutter $\cal C$ with $\cal D$ satisfying the above conditions yet.
These conditions will be used in Subsection \ref{subsec:projective}.

\section{Unique maximum fractional packing}\label{sec:unique}

In this section, we consider the problem under an additional condition that the maximum fractional packing is unique. Many important classes of precore clutters satisfy this condition.
Subsection \ref{subsec:unique} is concerned with step 1 in the case that the maximum fractional packing is unique. In Subsection \ref{subsec:projective}, we consider an example of a precore clutter. Moreover we show that there exists no counterexample to Conjecture \ref{C&G&M} in that class (Theorem \ref{thm:affinenon}).


\subsection{Unique maximum fractional packing}\label{subsec:unique}

In this subsection, we consider a clutter which has a unique maximum fractional packing.
For example, the clutter $Q_6=\{abc, cde, efa, bdf\}$ has a unique maximum fractional packing.

\begin{lem}\label{lem:unique}
Consider a clutter $\cal C$ which satisfies the tilde-full condition. The incidence vectors of $\tilde{\cal C}$ are affinely independent if and only if its maximum fractional packing is unique.
\end{lem}

\begin{pf}
Assume that a maximum fractional packing is unique.
Then $y(H)>0$ for all $H\in \tilde{\cal C}$ by the tilde-full condition.
The support of the maximum fractional packing $y$ of $\cal C$ consists of the hyperedges of $\tilde{\cal C}$ by Lemma \ref{lem:ibc}.  When these incidence vectors are affinely dependent, the maximum fractional packing can be moved slightly so that it is still a maximum fractional packing, a contradiction. 

When a maximum fractional packing is not unique, by taking two maximum fractional packings $y_1$ and $y_2$, they are affinely dependent since $\sum_{H\in {\cal C}}y_1(H) = \sum_{H\in {\cal C}}y_2(H)$ and $\sum_{H\in {\cal C}}y_1(H)1_H = \sum_{H\in {\cal C}}y_2(H)1_H=1_E$ by Lemma \ref{lem:mfp}. 
\end{pf}

Generally, when a polyhedron $P$ is not full dimensional, its facet-defining inequality is not unique. Here, we call a linear inequality $\langle 1_B,x \rangle \geq 0$ a {\it facet-defining inequality} of $P$ to a facet $F$ when $\{x\in P|\langle 1_B,x \rangle = 0\}=F$.

\begin{lem}\label{lem:simplex}
Consider a clutter $\cal C$ such that $\mbox{I}({\cal C})$ satisfies the integral blocking condition and is an integral polytope. 
Its maximum fractional packing is unique if and only if $\mbox{I}({\cal C})$ is a simplex.
\end{lem}

\begin{pf}
Since $\mbox{I}({\cal C})$ is a polytope, the minimum transversals of $\cal C$ cover $E$ by Lemma \ref{lem:nonempi}.
By Lemma \ref{lem:itos}, $\cal C$ satisfies the tilde-full condition.
By Lemma \ref{lem:unique}, its maximum fractional packing is unique if and only if the incidence vectors of $\tilde{\cal C}$ are affinely independent.
\end{pf}

We call an integral polytope which is simplex an {\it integral simplex}.
For $H\in \tilde{\cal C}$, we call a transversal $B\in \bl(\tilde{\cal C})$  a {\it facet transversal} of $H$ if $|H\cap B|>1$, $|B-H|\leq \bn({\cal C}) -2$, and $|H'\cap B|=1$ for $H'\in \tilde{\cal C}-\{H\}$.

\begin{thm}\label{thm:simplexfacet}
Consider the minimum-transversal-covered tilde-invariant clutter ${\cal C}$ which satisfies the integral blocking condition.
 $\mbox{I}({\cal C})$ is an integral simplex and ${\cal C}$ is hyperedge-nonseparable if and only if $\cal C$ satisfies the dimension condition and, for each $H\in {\cal C}$, there exists a facet transversal $B\in \bl({\cal C})$ of $H$.
\end{thm}

\begin{pf}
Assume that $\mbox{I}({\cal C})$ is an integral simplex. 
The extreme points of $\mbox{I}({\cal C})$ are the incidence vectors of the hyperedges of ${\cal C}$ by Lemma \ref{lem:tanten}. For each extreme point $x$ of $\mbox{I}({\cal C})$, we show that there exists a facet transversal $B\in \bl({\cal C})$ defining the facet which does not contain $x$. Since $\mbox{I}({\cal C})$ is an integral simplex, $x=1_H$ holds for some $H\in \tilde{\cal C}$.

Next we show that $|B-H|= \bn({\cal C}) -1$ does not hold. Assume on the contrary that $|B-H|= \bn({\cal C}) -1$.
Let $a\in B\cap H$. Then we have $\langle 1_{(B-H)\cup a}, 1_H \rangle =1$. Since  $\langle 1_{(B-H)\cup a},1_{H'} \rangle =1$ for $H'\in \tilde{\cal C}-\{H\}$,
$\langle 1_{(B-H)\cup a},x \rangle = 1$ for any $x\in \mbox{I}({\cal C})$.
Therefore $\langle 1_B,x \rangle \geq 1$ defines the same facet to $x(a)\geq 0$ since $\langle 1_{(B-H)\cup a},x \rangle = 1$. 

Assume that such a facet is defined by a linear inequality of the form $x(a)\geq 0$. The point $1_E$ is attained by the non-negative combination of ${\cal C}$ by Lemma \ref{lem:mfp}.
But every hyperedge in ${\cal C}$ except $H$ with $1_H=x$ does not contain $a$.
Therefore when $1_E$ is represented as a non-negative combination $y$ of $\cal C$, the coefficient $y(H)$ to $H$ is 1.
Since the deletion of all the elements reduces the fractional packing number by exactly one, it is hyperedge-separable.
Therefore any facet-defining inequality is defined by a facet transversal. 

The dimension condition follows from Lemma \ref{lem:dim}.

Conversely, assume that, for any $H\in {\cal C}$, there exists a transversal $B\in \bl({\cal C})$ such that $|H\cap B|>1$ and $|B-H| \leq \bn({\cal C})-2$ and $|H'\cap B|=1$ for $H'\in {\cal C}-\{H\}$. Then the incidence vector $1_H$ of every hyperedge $H$ in ${\cal C}=\tilde{\cal C}$ belongs to $\mbox{I}({\cal C})$ by Lemma \ref{lem:tanten}.
For each $1_H\in \mbox{I}({\cal C})$, the facet $\langle 1_B,x\rangle =1$ contains all the extreme points except $1_H$. By the dimension condition, the dimension of $\mbox{I}({\cal C})$ is the size $|{\cal C}|-1$. Therefore it becomes a simplex. Since its extreme points are expressed as ${\cal C}$, it is an integral simplex. 

Since $|B-H| \leq \bn({\cal C})-2$ holds, deleting all the elements of $H$ produces a clutter of blocking number at most $\bn({\cal C})-2$.
So ${\cal C}$ is hyperedge-nonseparable.
\end{pf}



We give an example of a precore clutter.
A graph $G$ is called a {\it brick} if it is 3-connected and $G-\{u,v\}$ has a perfect matching for all pairs of distinct $u,v\in V(G)$.
For a graph $G$, the {\it vertex cut clutter} ${\cal C}(G)$ is the clutter $\{\{\{a,b\}\in E(G)\} | a\in V(G)\}$.

\begin{ex}
For a brick $G$, the vertex cut clutter ${\cal C}(G)$ of $G$ is an example of a precore clutter.
In fact, since $G$ is non-bipartite, the maximum fractional packing is unique and the integral blocking condition is satisfied.
Since $G$ is matching covered, ${\cal C}(G)$ satisfies the minimum-transversal-covered.
Moreover, ${\cal C}(G)=\widetilde{{\cal C}(G)}$ holds.
Since the dimension of the matching polytope of a brick $G$ is $|E(G)|-|V(G)|$,
 the dimension condition is satisfied.
For any vertex $x\in V(G)$, there exists a factor of a brick $G$ such that vertex $x$ has degree 3 and the other vertices have degree 1. Such a factor becomes a facet-transversal. Therefore by Theorem \ref{thm:simplexfacet}, $I({\cal C}(G))$ is an integral simplex and hyperedge-nonseparable. We can show that ${\cal C}(G)$ is also non-separable by the definition of a brick.
\end{ex} 

\subsection{Combinatorial affine planes}\label{subsec:projective}

A clutter $\cal C$ on a finite set $E$ is called a {\it combinatorial projective plane} if the following three conditions are satisfied.

(1) For any two distinct elements, there exists a unique hyperedge containing the two elements.

(2) Any two distinct hyperedges intersect in exactly one element.

(3) There are four elements such that no line is incident with more than two of them.

On a combinatorial projective plane $\cal C$, each hyperedge is also called a {\it point} and each element is also called a {\it line}. (The inverse correspondence between point and line is possible but the reason why we adopt this correspondence is due to a combinatorial affine plane appeared later.)
For any combinatorial projective plane $\cal C$, there exists a natural number $n$ such that $n^2+n+1=|{\cal C}|=|E|$. Every element $a\in E$ is contained in $(n+1)$ hyperedges, and every hyperedge has size $n+1$.




By deleting one element $a\in E$ from the clutter $\cal C$, we obtain another clutter ${\cal C}\backslash a$.
This clutter ${\cal C}\backslash a$ becomes a clutter of a combinatorial affine plane. 

\begin{df}\normalfont
A clutter $\cal C$ on a finite set $E$ is called a {\it combinatorial affine plane} if the following three conditions are satisfied.

(1) For any two distinct hyperedges $H$ and $H'$, $|H\cap H'|=1$.

(2) Given an element $a$ and a hyperedge $H\in{\cal C}$ with $a\notin H$, there exists a unique element $b\in H$ such that $a$ and $b$ are not contained in the same hyperedge.

(3) There exist three hyperedges which do not contain the same element.
\end{df}

For example, a combinatorial projective plane on seven elements induces a combinatorial affine plane on six elements. It is $Q_6$, which is an ideal clutter of blocking number 2.


The next proposition is known.

\begin{prop}\label{prop:transversalhyperedge}
Let the size of a hyperedge of a combinatorial affine plane be $n+1$. Then $|E|=n^2+n$, the size of its minimum transversal is $n$, the number of its minimum transversals is $n+1$. Each element is contained in exactly $n$ hyperedges. Its minimum transversals form a partition of $E$. Any two elements which belong to different minimum transversals are included in exactly one hyperedge. Each hyperedge is also a transversal.
\end{prop}


\begin{lem}\label{lem:ti}
A combinatorial affine plane ${\cal C}$ satisfies the integral blocking condition, and ${\cal C}=\tilde{\cal C}$ holds.
Therefore it is tilde-invariant.
\end{lem}

\begin{pf}
Let $n+1$ be the size of its hyperedge.
We first show the integral blocking condition. The blocking number of $\cal C$ is $n$. For each element in $E$, there exist $n$ hyperedges of ${\cal C}$ containing the element. Therefore the sum of the incidence vectors of all the hyperedges of ${\cal C}$ is $n1_E$, which is a fractional packing of ${\cal C}$. Since they form a maximum fractional packing, ${\cal C}$ satisfies the integral blocking condition.

Since any minimum transversal and any hyperedge of ${\cal C}$ intersect in exactly one element by Proposition \ref{prop:transversalhyperedge}, we have ${\cal C}=\tilde{\cal C}$. 
\end{pf}

\begin{lem}\label{lem:punique}
The maximum fractional packing of a combinatorial affine plane is unique.
\end{lem}

\begin{pf}
By calculating the determinant of the clutter matrix, the incidence vectors of hyperedges of the clutter are affinely independent.
So the statement follows from Lemma \ref{lem:unique}.
\end{pf}

\begin{thm}\label{thm:cppintegral}
For a combinatorial affine plane ${\cal C}$, $\mbox{I}({\cal C})$ is an integral simplex and ${\cal C}$ is nonseparable.
Therefore $\cal C$ is a precore clutter.
\end{thm}

\begin{pf}
We can take the hyperedges on $\cal C$ as facet transversals by Proposition \ref{prop:transversalhyperedge}.
The number of the hyperedges is $n^2$ and the number of the minimum transversals is $n+1$.
Since they are affinely independent, we have the dimension condition.
By Theorem \ref{thm:simplexfacet}, $\mbox{I}({\cal C})$ is an integral simplex.

By deleting all the points on a hyperedge, all the hyperedges disappear. So the clutter of a combinatorial affine plane is  nonseparable.
\end{pf}


\begin{thm}\label{thm:affinenon}
Every combinatorial affine plane ${\cal C}$ of blocking number at least 3 has no ideal minimally non-packing solution clutter.
\end{thm}

\begin{pf}
Assume that ${\cal C}$ has an ideal minimally non-packing solution clutter $\cal D$.

Consider distinct hyperedges $A,B, $ and $C$ in ${\cal C}$ with $A\cap B\cap C=\emptyset$. Let $z$ be a unique point in $A\cap B$. Similarly, let $x$ be a unique point in $B\cap C$, and $y$ be a unique point in $C\cap A$. 

\begin{figure}[ht]
\begin{center}
\includegraphics[width=4cm]{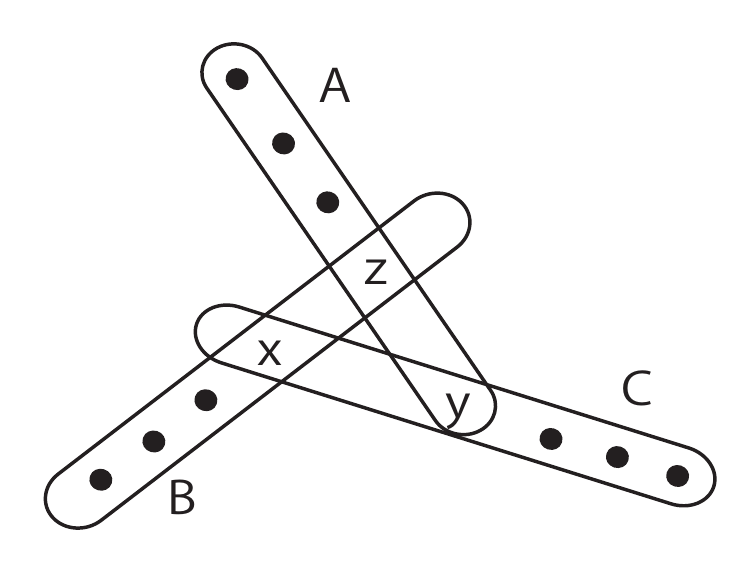}
\caption{an example of ${\cal D}[X]$}
\end{center}
\end{figure}

Then consider the restriction ${\cal C}[X]$ where $X$ is the union of the three hyperedges $A, B,$ and $C$. Note that $X\neq E$ since the blocking number is at least 3. Then since such a clutter has exactly three hyperedges, its blocking number is 2. Since ${\cal D}[X]$ must pack, there exists a packing of size 2 in ${\cal D}[X]$ (Condition B). By Condition IF, any facet transversal of $I({\cal C})$ is also a facet transversal of $I({\cal D})$. Therefore any hyperedge $H\in {\cal C}$ is also a transversal in ${\cal D}$. Moreover any minimum transversal of ${\cal C}$ is a minimum transversal of ${\cal D}$ by Condition IM. Therefore the two elements consisting of $x$ and any one element of $A-\{y,z\}$ form a transversal of ${\cal D}[X]$. Similarly,  two elements consisting of $y$ and any one element of $B-\{z,x\}$ form a transversal, and two elements consisting of $z$ and any one element of $C-\{x,y\}$ form a transversal. Such two elements are included in some minimum transversal or included in some hyperedge which is also a transversal in $\bl({\cal D})$ since the deletion of elements from a clutter corresponds to the contraction of them from the clutter of transversals. Therefore $X$ is covered by transversals of size 2. By regarding such two elements as an edge of a graph, such a graph has three connected components and each of them is a star. A packing of size 2 becomes a partition on $X$ consisting of two hyperedges of size $|X|/2$ as in Figure \ref{fig:stars}. For a packing of size 2 in ${\cal D}[X]$, two elements as a transversal belong to different hyperedges in the packing of size 2 on ${\cal D}[X]$. Therefore we can take four types of packings of size 2.

\begin{figure}[ht]
\begin{center}
\includegraphics[width=6cm]{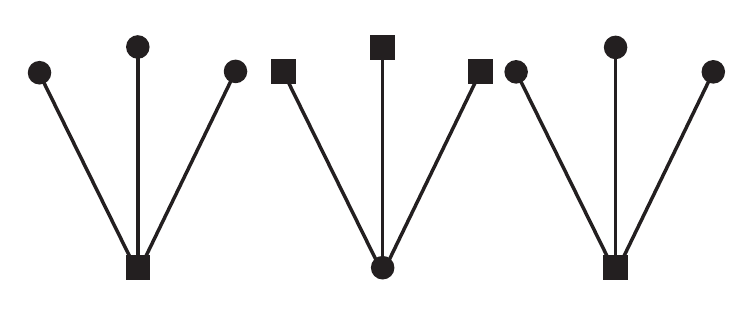}
\caption{sets of vertices indicated by circles and squares represent hyperedges}
\label{fig:stars}
\end{center}
\end{figure}

In three types out of the four types of packings, one hyperedge in packings of size 2 is either of $A, B, $ and $C$, the other hyperedge is included in the complement of the hyperedge in $X$.
These cases contradict the fact that $A,B,C$ themselves are transversals because any hyperedge must intersect any transversal. We discuss the remaining type of the packings, that is, one hyperedge is included in $\{x,y,z\}$, and other hyperedge is included in $X-\{x,y,z\}$. Since $\{x,y,z\}$ intersects any minimum transversal in exactly one element, $\{x,y,z\}$ cannot be a hyperedge of $\cal D$ by Condition H, a contradiction.
\end{pf}

We should note that whether 
a combinatorial projective plane except for the $F_7$ can 
be a core of a minimally non-ideal clutter or not is a famous 
open question of the theory of minimally non-ideal clutters 
(see Question 9 in Cornu\'{e}jols, Guenin and Tun\c{c}el \cite{Cornuejols_Guenin_and_Tuncel_2008}).
Conjecture \ref{con:affineideal} implies their conjecture.

\begin{con}\label{con:affineideal}
A clutter of a combinatorial affine plane of blocking number at least 3 has no ideal solution clutter.
\end{con}






















\end{document}